\documentclass{article}
\usepackage{amsmath,amssymb, amscd,epsfig,latexsym,amsfonts, ebezier,authblk}
\usepackage{amsmath,amsthm,verbatim,amssymb,amsfonts,amscd, graphicx}
\usepackage{amsfonts,color,xcolor,framed,amsthm}
\usepackage{sidecap}
\usepackage{booktabs} 
\usepackage{graphicx}
\usepackage{blindtext}
\usepackage{pdfsync}
\usepackage{subfigure,wrapfig,sidecap}
\usepackage{color}
\usepackage{epstopdf}
\usepackage{algorithmic}

\def\d{{\rm d}}

\newcommand{\be}{\begin{equation}}
\newcommand{\eeq}{\end{equation}}
\newcommand{\bea}{\begin{eqnarray}}
\newcommand{\eea}{\end{eqnarray}}

\newtheorem{conjecture}{Conjecture}
\newtheorem*{remark}{Remark}

\DeclareMathOperator{\sech}{sech}
\DeclareMathOperator{\ee}{e}
\DeclareMathOperator{\ii}{i}
\DeclareMathOperator{\de}{d}

\title{ Numerical study of the Kadomtsev--Petviashvili equation and  dispersive shock waves}
\author[1,2]{T. Grava}
\author[3]{C. Klein}
\author[1]{G. Pitton}
\affil[1]{Scuola Internazionale Superiore di Studi Avanzati, Trieste, Italy}
\affil[2]{School of Mathematics, University of Bristol, UK}
\affil[3]{Institut de Math\'ematiques de Bourgogne, Universit\'e de 
Bourgogne-Franche-Comt\'e, France}
\begin{document}
 \maketitle

%

\begin{abstract}
  A detailed numerical study of the long time behaviour of 
    dispersive shock waves in solutions to the Kadomtsev-Petviashvili 
    (KP)
    I equation  is presented. It is 
    shown that   modulated lump solutions emerge from  the   dispersive shock waves.  For the description of dispersive shock waves,  Whitham modulation equations for KP  are obtained. It is 
    shown that the modulation  equations  near the   soliton line   are hyperbolic  for the KPII   equation while they are elliptic  for the KPI equation 
    leading  to a focusing effect and the formation of lumps. Such a 
    behaviour is similar  to the appearance of breathers for   the focusing 
    nonlinear Schr\"odinger equation in the semiclassical limit. \end{abstract}

\section{Introduction}
\begin{figure}[h]
  \centering
   \includegraphics[width=0.45\textwidth]{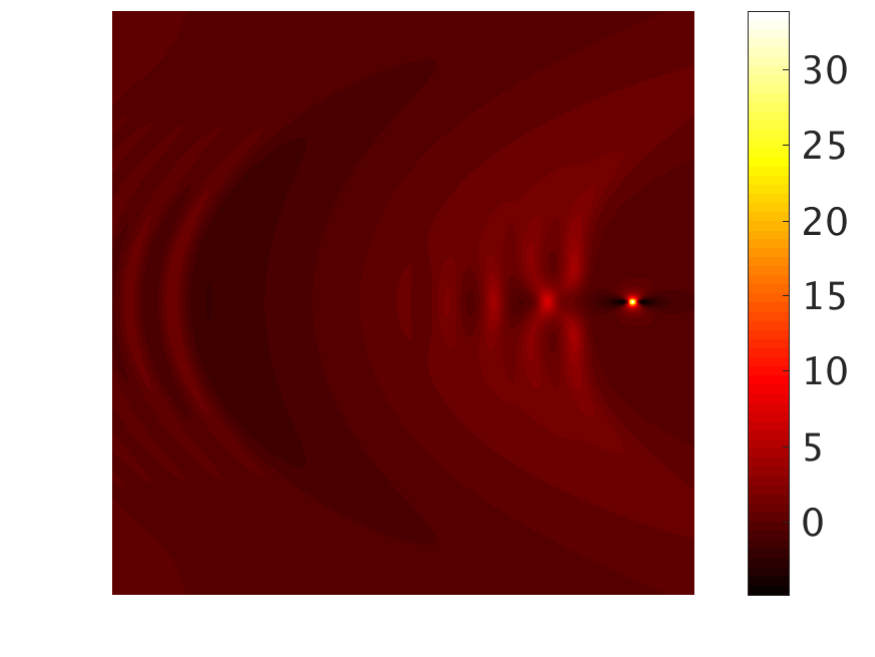}
  \includegraphics[width=0.45\textwidth]{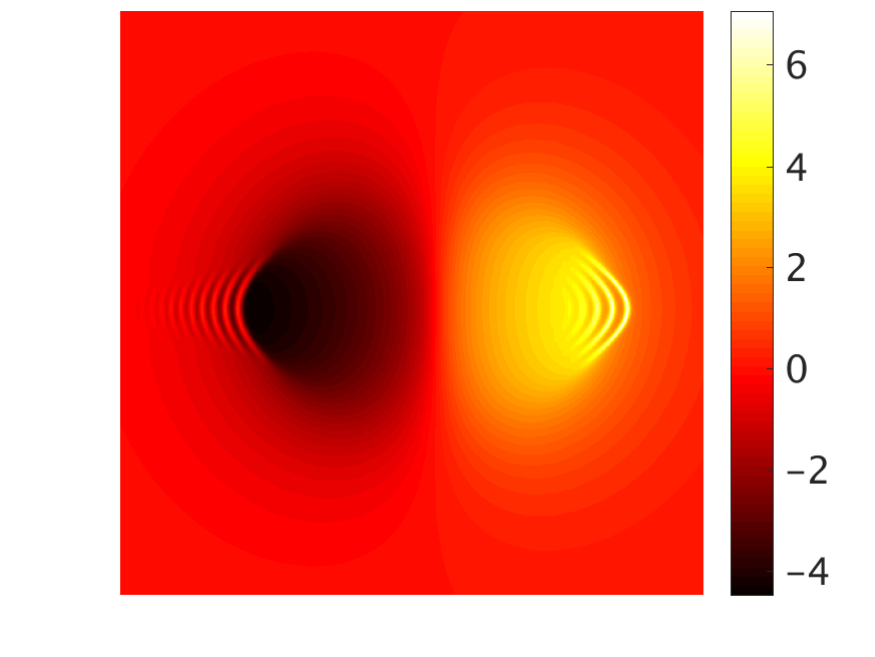} \\
  \includegraphics[width=0.49\textwidth]{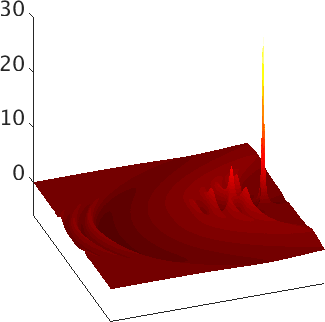}
  \includegraphics[width=0.49\textwidth]{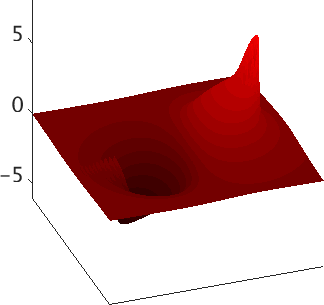}
  \caption{Solution of the KPI equation (left) and of the KPII 
  equation (right) for $\epsilon=0.1$ and the initial data 
  $u(x,y,0)=-6\partial_{x}\mbox{sech}^{2}x$ at time $t=0.8$. Notice how the 
  KPI solution has developed 
  a spike that is about 5 times higher than the highest peak of the  KPII solution.}
  \label{KPI_KPII}
\end{figure}

We consider the  Cauchy problem for  the Kadomtsev  Petviashvili  (KP) equation
\begin{equation}
\label{KP}
(u_t+uu_x+\epsilon^2 u_{xxx})_{x}+\alpha u_{yy}=0,\quad \alpha=\pm 1,
\end{equation}
in the class of rapidly decreasing smooth initial data.
Here  $\epsilon>0$ is a small parameter  and we are interested in the behaviour of the solution $u(x,y,t;\epsilon)$ as $\epsilon\to 0$.  In such a limit the solution of the KP equation develops strong
oscillations and very high peaks that will be the subject of the present manuscript.
 The equation (\ref{KP})   was  first introduced by Kadomtsev and Petviashvili  
 \cite{KP70} in order to study the stability of the Korteweg--de Vries 
 (KdV)
 soliton in a two-dimensional setting, and  it  is now a prototype  for the evolution of  weakly nonlinear
  quasi-unidirectional  waves of small amplitude  in various physical 
  situations.  
  For $\alpha=-1$  ($\alpha=1) $ the  equation (\ref{KP}) is called KPI (KPII)  equation and  describes   quasi-unidirectional long waves in shallow water with weak transversal effects and 
 strong  (weak) surface tension. 
  The KPII equation is known to have a defocusing effect, 
  whereas the KPI equation is focusing. It is exactly this latter 
  effect which we will study in this paper. A comparison of the 
  solutions of the two KP equations for the same initial data is 
  shown in Fig.~\ref{KPI_KPII} where one can see the focusing effect of KPI.

  The  KP equation is also the prototypical  integrable equation \cite{Dryuma} in two spatial dimensions  and it has been studied via inverse scattering  \cite{FokasAblowitz} \cite{BP}.
In the dimensionless  KP equation, i.e., equation (\ref{KP}) with 
$\epsilon=1$,  a
parameter $\epsilon$  is  introduced by considering the  long time behavior of solutions with  slowly varying initial data of the form  $ u_0(\epsilon x,\epsilon y)$  where $0 < \epsilon\ll 1$
is a 
small parameter and $ u_0(x,y)$  is some given initial profile. As  $\epsilon \to 0$  the initial datum  approaches a constant   value
and  in order to see nontrivial effects one has to wait for times  of order $ t \simeq  O(1/\epsilon),$ 
 which consequently requires to rescale the spatial variables onto macroscopically large scales $ x \simeq O(1/\epsilon)$, too. 
 This is equivalent to  consider the rescaled variables  $x \to x' = x\epsilon,$ $y\to y' = y\epsilon $,  $t  \to t' = t\epsilon$  and put  $u^{\epsilon}(x',y',t') = u(x\epsilon,y\epsilon, t\epsilon)$
to obtain the equation (\ref{KP}) where we omit the $'$ for simplicity.

For $\epsilon=0$ the KP equation turns into the so called dispersion-less KP   equation (dKP)  \cite{LRT} ,\cite{ZK69} 
\begin{equation}
\label{dKP}
(u_t+uu_x)_x+\alpha u_{yy}=0.
\end{equation}
Note that in spite of its name, the dKP equation (\ref{dKP}) contains dispersion, and only the highest order dispersive term has been dropped relative to (\ref{KP}). 
Local well-posedness of the Cauchy problem for the dKP equation  has been proved in certain Sobolev spaces in \cite{Rozanova}.
Generically, the solution of the dKP equation develops a singularity 
in finite time  $t_c>0$.  It is discussed 
in \cite{GKE} and \cite{MS08}  that this singularity develops in a point where the gradients become divergent in all directions except one.

As long as the gradients of the dKP solution remain bounded, the solution $u(x,y,t;\epsilon)$ of the KP equation is expected to be approximated in the limit $\epsilon\to 0$ by the solution of the dKP equation.  Even if there are  many strong results about the Cauchy problem for the KP 
equation in various functional spaces  (see, e.g., 
\cite{Bourgain,  MST}), these results are  insufficient to rigorously justify the small $\epsilon$ behaviour of solutions to KP even for $0<t<t_c$.
Near $t=t_c$ the solution of the KP equation,  preventing  the 
formation of the strong gradients in the dKP solution, starts to develop a region of rapid modulated oscillations.
These oscillations are called dispersive shock waves, and they can be approximated  at the onset of their formation by a particular solution of the Painlev\'e I2 equation, up to shifts and rescalings  \cite{DGK}.

For later times $t>t_c$  these oscillations are expected to be described by the modulated  travelling  cnoidal  wave solution of the KP equation.
The travelling cnoidal wave solution is given by 
\begin{equation}
\label{cnoidal0}
u(x,y,t;\epsilon)=\beta_1+\beta_3-\beta_2+2(\beta_2-\beta_3)\mbox{cn}^2\left(\dfrac{\sqrt{\beta_1-\beta_3}}{\sqrt{6}\epsilon}( x+\frac{l}{k}y-\frac{\omega}{k} t)+\phi_0; m\right)
\end{equation}
where $\mbox{cn}(z; m)$ is the Jacobi elliptic function of modulus 
$m=\dfrac{\beta_2-\beta_3}{\beta_1-\beta_3}$ with the constants $\beta_1>\beta_2>\beta_3$,   $\phi_0$ is an arbitrary constant and $K(m)$ the complete elliptic integral of the first kind.
The wave number $k$  and the frequency $\omega $ are given by  
\begin{equation}
\label{wavenumber}
k=\pi \dfrac{\sqrt{\beta_1-\beta_3}}{\sqrt{6}K(m)},\quad 
\omega=\dfrac{k}{3}(\beta_1+\beta_2+\beta_3)+\alpha\dfrac{l^2}{k}.
\end{equation}
The average value $\bar{u}$ over a period and the maximum amplitude $a:=u_{max}-u_{min}$ of the oscillations are 
\begin{equation}
\label{u_ave}
\bar{u}=\beta_2+\beta_3-\beta_1+2(\beta_1-\beta_3)\dfrac{E(m)}{K(m)},\quad a=2(\beta_2-\beta_3),
\end{equation}
where   $E(m)$  is the complete elliptic integral of the second kind.
For constant values of $\beta_1,\beta_2,\beta_3$ and $l$, the formula (\ref{cnoidal0})  gives an exact solution of the KP equation.
%
The modulation of the wave-parameters  of the cnoidal wave solution is obtained by letting $\beta_j=\beta_j(x,y,t)$, $j=1,2,3$  and $l=l(x,y,t)$ and   requesting that (\ref{cnoidal0}) is an approximate solution of KP up to higher order corrections.
Over the last forty years, since the seminal paper of Gurevich and Pitaevsky,  \cite{GP}  there has been a lot of attention to the quantitative study of dispersive shock waves  see e.g. the recent volume
\cite{miller}, and refined experiments have been developed \cite{trillo}.
Most of the analysis is restricted 
to models in one spatial dimension. Two dimensional models have been 
much less studied, see for example \cite{Hoefer},\cite{Ratliff},\cite{El}.  Regarding the KP 
equation,  the formation of dispersive shock waves has been studied 
numerically in \cite{KSM,KR13} and  both numerically and analytically  in  \cite{Ablowitz}  for  an 
initial step with parabolic profile, and recently in \cite{Biondini} using the method of multiple scales. Modulation theory in the general setting of Riemann surfaces has been developed 
in \cite{Krichever}.
In this manuscript we derive  the modulation equations for  KP  using the Whitham averaging method over the Lagrangian as in \cite{Whitham}.
Our final form of the equations  for $\beta_1(x,y,t)>\beta_2(x,y,t)>\beta_3(x,y,t)$ and $q(x,y,t):=l(x,y,t)/k(x,y,t) $, plus two extra dependent variables $p=p(x,y,t)$ and $r=r(x,y, t)$   (see definition (\ref{q} and (\ref{def_r}))  is 
\begin{align}
\label{eqb_in}
&\dfrac{\partial}{\partial t}\beta_i+(v_i+\alpha q^2)\dfrac{\partial }{\partial x}\beta_i+\alpha(2 q  D\beta_i-(v_i-2\beta_i)Dq+Dp)=0,\quad i=1,2,3,\\
\label{eqq1_in}
&\dfrac{\partial}{\partial t}q+\left(
\dfrac{1}{3}\sum_{i=1}^3\beta_i+\alpha q^2\right)q_x+2\alpha 
Dq+\dfrac{1}{3}D(\sum_{i=1}^3\beta_i)=0,\\
\label{eep_in}
&p_t+\left(
\dfrac{1}{3}\sum_{i=1}^3\beta_i+\alpha q^2\right)p_x+Dr=0,\;\;
r_x-\dfrac{B_x}{6}-\alpha(\bar{u}(q_y-qq_x)+p_y+qp_x))=0,
\end{align}
with $D=\dfrac{\partial}{\partial y}-q\dfrac{\partial}{\partial x}$, the speeds $v_i=v_i(\beta_1,\beta_2,\beta_3)$ are 
\begin{equation}
\label{vi}
v_i=\dfrac{1}{3}(\beta_1+\beta_2+\beta_3)+\dfrac{2}{3}\dfrac{\prod_{k\neq i}(\beta_i-\beta_k)}{\beta_i-\beta_1+(\beta_1-\beta_3)\dfrac{E(m)}{K(m)}}, \quad i=1,2,3,
\end{equation}
with  $\bar{u}$ defined in (\ref{u_ave}) and  
$B=\sum_{i=1}^3\beta_i^2-2(\beta_1\beta_2+\beta_2\beta_3+\beta_1\beta_3)$.  The system satisfies two compatibility conditions given by the constraints
\begin{equation}
\label{const_final_in}
q_x=\dfrac{k_y}{k}-q\dfrac{k_x}{k},\;\;\;p_x=\bar{u}_{y}-(q\bar{u})_x.\;\;\;
\end{equation}
When $\alpha=0$ the equations (\ref{eqb_in})   and the second 
equation in (\ref{eep_in}) coincide with the Whitham modulation 
equations for KdV with $r=B/6$ being an integral.
The equations (\ref{eqb_in}), (\ref{eqq1_in}) and  (\ref{const_final_in})  are equivalent to  the equations obtained in \cite{Biondini},  while the equations (\ref{eep_in}) seem to be new. We set up the Cauchy problem for the Whitham modulation equations and we show that the Whitham system near the solitonic front   when $m\simeq 1$ is not hyperbolic.

When the modulus $m\to 1$, the travelling wave solution  (\ref{cnoidal0})  of KP converges to 
\begin{equation}
\label{soliton}
u(x,y,t;\epsilon)\simeq \beta_3+2(\beta_1-\beta_3)\mbox{sech}^2\left(\frac{\sqrt{\beta_1-\beta_3}}{\sqrt{6}\epsilon}(x+\frac{l}{k}y-\frac{\omega}{k} )t+\phi_0\right).
\end{equation}
If we set $\beta_3=0$ and $\beta_1=6k^2$, the above expression is exactly the line soliton of the KP equation and the wave numbers $k$, $l$ and $\omega$ satisfy the  dispersion  relation
$\omega=4k^3+\dfrac{\alpha l^2}{k}$ (see e.g. \cite{AC}).
For the KPI equation the line soliton is known to be linearly 
unstable under perturbations,  \cite{Zak}, \cite{RT}.
Numerical studies as \cite{Infeld92},  see also the more recent papers  \cite{KS12,KSM},  and analytical studies \cite{Pelinovsky} indicate that the  solitons of the form (\ref{soliton}) of 
sufficient amplitude are unstable against the formation of so called  lump solutions.

Lumps are  localised solutions  decreasing algebraically at infinity that  take the form
 \begin{equation}
 \label{lumpgeneral}
 u(x,y,t;\epsilon)=
 24\dfrac{(-\frac{1}{\epsilon^2}(x+ay +(a^2-3b^2)t)^2+3\frac{b^2}{\epsilon^2}(y+2at)^2+1/b^2)}{(\frac{1}{\epsilon^2} (x+ay +
 (a^2-3b^2)t)^2+3\frac{b^2}{\epsilon^2}(y+2a t)^2+1/b^2)^2},
 \end{equation}
 where  $a$ and  $b$  are  arbitrary constants.
The maximum of the lump is located at 
  \[
 x=3b^2t+a^2t,\quad y=-2at,  \]
 with maximum value $24b^2$. When $a=0$ the lump is symmetric with respect to $y$-axis.
%

\begin{figure}[htbp]
  \centering
 
   \subfigure[$\epsilon=0.05$, $C_0=6$, $t=0.8$]{
  \includegraphics[width=0.45\textwidth]{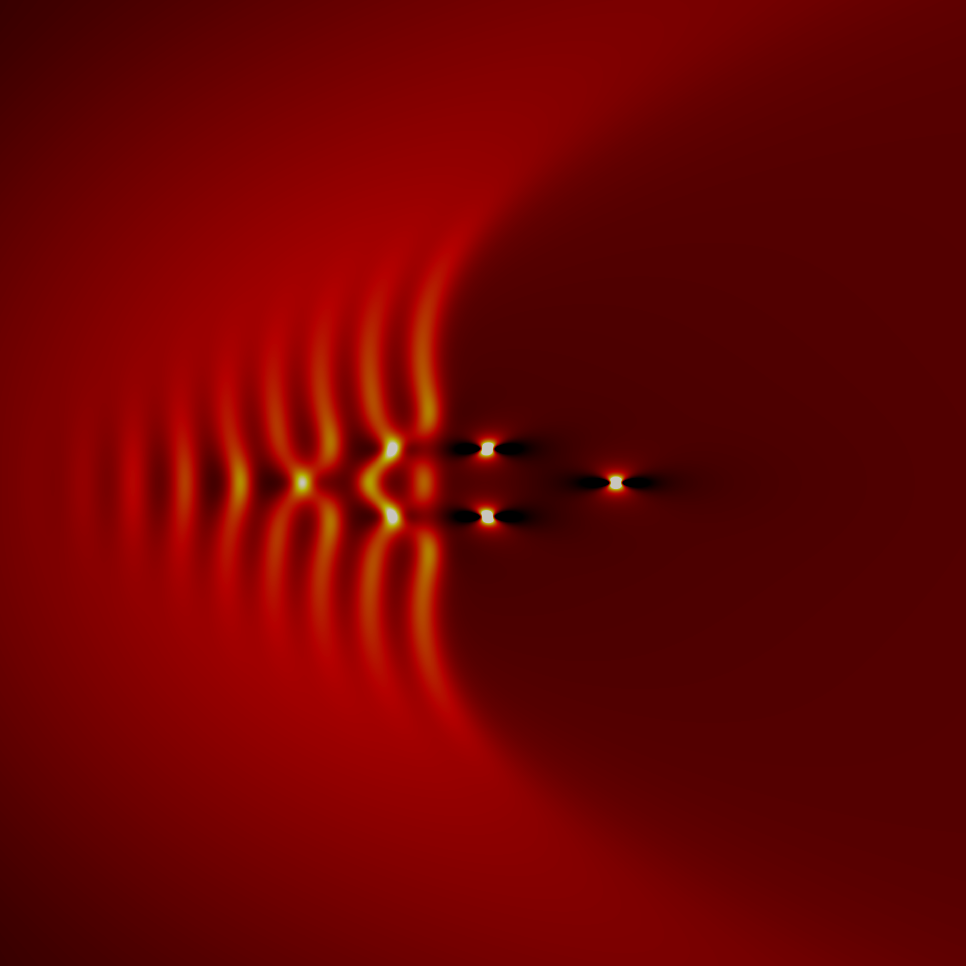}
   }
   \subfigure[$\epsilon=0.02$, $C_0=6$, $t=0.8$]{
      \includegraphics[width=0.45\textwidth]{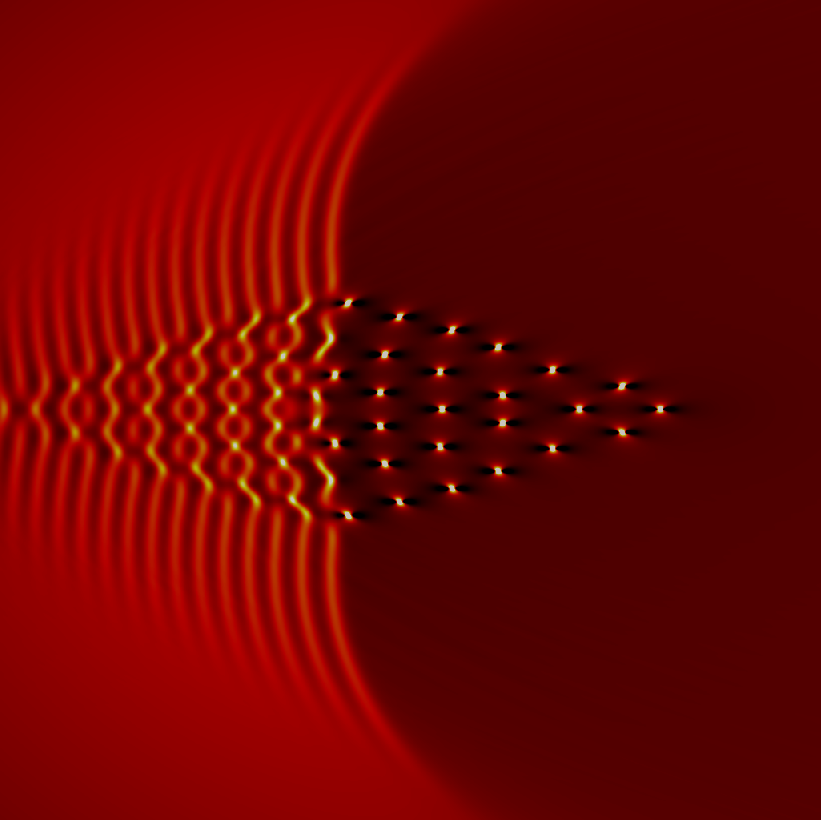}
  }
  \caption{Detail of the lattice arrangement of the lumps  on the $(x,y)$ plane for two representative cases  for the initial data $u_0(x,u)=-C_0\partial_x\sech^2\sqrt{x^2+y^2}$. The distance between the lumps clearly scales with $\epsilon$.}
  \label{soliton_lattice}
\end{figure}

We obtain, using the averaging over Lagrangian density,  the 
modulation of the soliton parameters. These equations are elliptic  for KPI and therefore they are expected 
to develop a point of elliptic umbilic catastrophe   as 
for the semiclassical limit of the focusing nonlinear Schr\"odinger 
(NLS) equation   \cite{DGK}. In the NLS case 
 a train of Peregrine breathers is  generically  formed  
\cite{Bertola_Tovbis}  that is  in amplitude three times the value  of the solution at the point of elliptic umbilic catastrophe.   Furthermore the position of the breathers scales in $\epsilon$ with the power $4/5$.
The soliton front of the dispersive shock waves for KPI  breaks 
 into a lattice of lumps and the distance among the lumps scales with $\epsilon$, see Fig.~\ref{soliton_lattice}

 The amplitude of the first lump that appears is proportional to the initial data and,  for the specific initial data considered,  it is  about ten times the maximal amplitude of the initial data.
 The amplitude of the lump decreases (numerically)  with time,  without producing any radiation as in \cite{Minzoni}.
 Finally we study the dependence on $\epsilon$ of the position and the time of formation of the first lump and we find a scaling exponent that is compatible with the value $4/5$ as in the NLS case.

 This manuscript is organised as follows.
 In section 2 we derive the Whitham modulation equations for KP using 
 the averaging over the Lagrangian.
 We then define the 
 Cauchy problem for the Whitham modulation equations.
 Next we obtain the modulation equations of the soliton parameters 
 and show that for KPI such equations are elliptic. We then show that  the Whitham system is not hyperbolic near the soliton front, since two eigenvalues  of the velocity matrix are complex.
 In section 3 we  collect known results on the focusing NLS equation and on how 
 solutions to the NLS equation are related to KPI solutions. In 
 section 4 we briefly present the numerical methods used  for the 
 integration of the KP equation. These methods are applied in section 5 to 
 concrete examples for the KPI equation. In particular we study 
 numerically the nature of the lattice of  lumps that is formed out 
 of the soliton front  in the KPI solution in the small dispersion 
 limit. We add some concluding remarks in section 6. 

\section{Whitham modulation equations for KP via Lagrangian averaging}
In this section we will obtain the Whitham modulation 
equations for the KP equations   following Whitham   method \cite{Whitham} of averaged Lagrangian as in   \cite{Infeld92}.

\subsection{ Lagrangian  density for the travelling wave solution of KP}
The Lagrangian density of the KP 
equation is 
\begin{equation}
\label{Lagrangian}
L=\epsilon^2f_tf_x+\dfrac{\epsilon^3}{3}f_x^3-\epsilon^4f_{xx}^2+\epsilon^2\alpha f_y^2
\end{equation}
which leads to the Euler-Lagrange equation
\[
f_{tx}+\epsilon f_x f_{xx}+\epsilon^2 f_{xxxx}+\alpha f_{yy}=0.
\]
The above equation  coincides with the KP equation under the substitution
$\epsilon f_x=u$.
We look for a solution that is a travelling wave, namely a  solution of the form
\[
f=\psi+\phi(\theta),\quad \theta=\frac{kx+ly -\omega  t}{\epsilon},\quad \psi=\frac{c_1x+c_2 y-\gamma t}{\epsilon},
\]
where  $\phi(\theta)$ is  a $2\pi$ periodic function of its argument and the remaining quantities are parameters to be determined.  In our notation $x/\epsilon$, $y/\epsilon $ and $t/\epsilon$ will be the fast variables and 
$x,y$ and $t$ will be the slow variables. We introduce
\[
\eta=\epsilon f_x=c_1+k\phi_{\theta},
\quad 
\epsilon f_y=c_2+\dfrac{l}{k}(\eta-c_1),\quad 
\epsilon f_t=-\gamma-\dfrac{\omega}{k}(\eta-c_1).
\]
It follows from (\ref{KP}) that the function $\eta(\theta)$ satisfies the equation
\begin{equation}
\label{spectralKP}
3k^2\eta_{\theta}^2=-\eta^3+V\eta^2+B\eta+A,\
\end{equation}
where  $B$ and $A$ are integration constants and 
\begin{equation}
\label{V}
V=3\left(\dfrac{\omega}{k}-\alpha\dfrac{l^2}{k^2}\right).
\end{equation}
 In order to get a periodic solution, we assume that the polynomial 
 \begin{equation}
 \label{ei}
 -\eta^3+V\eta^2+B\eta+A=-(\eta-e_1)(\eta-e_2)(\eta-e_3)
 \end{equation}
with $e_1>e_2>e_3$. Then the periodic motion takes place for $e_2\leq \eta\leq e_1$ and one has the relation
\begin{equation}
\label{spectral1}
\sqrt{3}k\dfrac{\de\eta}{\sqrt{(e_1-\eta)(\eta-e_2)(\eta-e_3)}}=\de\theta,
\end{equation}
so that integrating over a period, one obtains 
\[
2\sqrt{3}k\int^{e_1}_{e_2}\dfrac{\de\eta}{\sqrt{(e_1-\eta)(\eta-e_2)(\eta-e_3)}}=\oint \de\theta=2\pi.
\]
It follows  that the wave number $k$ can be expressed in terms of  a complete integral of the first kind:
\begin{equation}
\label{k}
k=\pi\dfrac{\sqrt{(e_1-e_3)}}{2\sqrt{3}K(m)},\quad m=\dfrac{e_1-e_2}{e_1-e_3},\quad K(m):=\int_0^{\frac{\pi}{2}}\dfrac{\de\psi}{\sqrt{1-m^2\sin^2\psi}}.
\end{equation}
Integrating between $e_2$ and $\eta$  in equation (\ref{spectral1}) 
one arrives to  the  expression
\begin{equation}
\label{cnoidal}
u(x,y,t)=\eta(\theta)=e_2+(e_1-e_2)\mbox{cn}^2\left(\frac{\sqrt{e_1-e_3}}{2\sqrt{3}\epsilon}\left(x-\frac{\omega}{k}t+\frac{l}{k} y\right)-K(m);m\right),
\end{equation}
where we use also the evenness of the function $\mbox{cn}(z;m)$.
The Lagrangian corresponding to the traveling wave solution  (\ref{cnoidal}) derived above  takes the form
\begin{equation}
\begin{split}
\label{L}
L=-2k^2\eta_{\theta}^2+\eta\left(\frac{B}{3}-\gamma+c_1\dfrac{\omega}{k}+2\alpha\dfrac{l}{k}\left(c_2-\dfrac{l}{k}c_1\right)\right)+\alpha\left(c_2-\dfrac{l}{k}c_1\right)^2+\frac{A}{3}.
\end{split}
\end{equation}

\subsection{Whitham average equations via Lagrangian averaging}
Below  we are going to apply Whitham's procedure to obtain the modulation of the wave parameters $A$, $B$, $V$, $k$, $l$, $c_1$,  $c_2$ and $\gamma$
by  variation of averaged quantities.
We introduce the averaged quantities
\begin{equation}
\label{average}
\langle \eta\rangle=\dfrac{1}{2\pi}\int_0^{2\pi}\eta \de\theta=c_1,\quad 
\langle \eta_{\theta}^2\rangle=\dfrac{1}{2\pi}\int_{0}^{2\pi} \eta_{\theta}^2\de \theta=\dfrac{W}{k},
\end{equation}
where 
 \[
 W:=\frac{1}{\sqrt{3}\pi}\int^{e_1}_{e_2}\sqrt{-\eta^3+V\eta^2+B\eta+A}\de\eta.
 \]
 Using (\ref{average}),  the average  of the  Lagrangian  $L$ defined in (\ref{L})   takes the form 
 \[
 {\mathcal L}:=
 \dfrac{1}{2\pi}\int_{0}^{2\pi} L\de\theta=-2kW+\dfrac{1}{3}Bc_1-\gamma c_1+\dfrac{1}{3}V c_1^2+\alpha c_2^2+\dfrac{A}{3}.
\]
The Lagrangian $ {\mathcal L}= {\mathcal L}(\omega,k,l,A,\gamma, c_1,c_2, B)$ and the Whitham method consists in assuming that the quantities 
$\omega,k,l,A,\gamma, c_1,c_2$ and $B$ 
 depend on the slow variables $ x$, $ y$ and $ t $. The variational principle is
 \[
 \delta\int\int  {\mathcal L}(\omega,k,l,A,\gamma, c_1,c_2, B)\de x\de y \de t=0.
 \]
 The variational equations are 
 (see (14.69)-(14.73) in \cite{Whitham})
 \begin{align}
\label{k_W}
&{\mathcal L}_A=0\to kW_A=\dfrac{1}{6},\quad {\mathcal L}_B=0\to \dfrac{ c_1}{6}=kW_B,\\
\label{W2}
&\dfrac{\partial }{\partial t} {\mathcal L}_{\omega}-\dfrac{\partial }{\partial x} {\mathcal L}_k-\dfrac{\partial }{\partial y} {\mathcal L}_l=0,\quad \dfrac{\partial }{\partial t} {\mathcal L}_{\gamma}-\dfrac{\partial }{\partial x} {\mathcal L}_{c_1}-\dfrac{\partial }{\partial y} {\mathcal L}_{c_2}=0
\end{align}
together with the consistency conditions  which follows from $\theta_{xt}=\theta_{tx}$, $\psi_{xt}=\psi_{tx}$, $\theta_{xy}=\theta_{yx}$, $\psi_{xy}=\psi_{yx}$ and 
$\theta_{yt}=\theta_{ty}$, $\psi_{yt}=\psi_{ty}$
\begin{align}
\label{KP1}
&k_t+\omega_x=0,\quad \dfrac{\partial}{\partial t} c_1+\dfrac{\partial}{\partial x}\gamma=0, \\
\label{KP2}
&l_t+\omega_y=0,\quad  \dfrac{\partial}{\partial t}c_2+\dfrac{\partial}{\partial y}\gamma=0,\\
\label{KP3}
&l_x=k_y,\quad\quad \dfrac{\partial}{\partial x} c_2= \dfrac{\partial}{\partial y}c_1.
\end{align}
Since KP can be written in the form
\begin{align}
&u_t+uu_x+\epsilon^2 u_{xxx}+\alpha v_y=0,\quad v_x=u_y,
\end{align}
one has, for the travelling wave
$kv_{\theta}=lu_{\theta}$,
which after integration in $\theta$ gives
$kv(\theta)=lu(\theta)+c_0$
for some integration constant $c_0=c_0(x,y,t)$ independent from $\theta$.
 Therefore we define the new dependent variables $p=p(x,y,t)$ and $q=q(x,y,t)$ as
\begin{equation}
\label{q}
 q:=\dfrac{l}{k},\quad p:=\langle v\rangle-q\langle u\rangle=c_2-qc_1\,.
\end{equation}
To simplify further the final form of the equations we also introduce  a new dependent variable  $r$  in place of $\gamma$
\begin{equation}
\label{def_r}
r:=\gamma-\dfrac{\omega}{k}c_1\,.
\end{equation}
Using (\ref{V}),  (\ref{k_W}) and (\ref{KP3}),    and the above definitions  we can write  the six equations (\ref{W2}),  (\ref{KP1}),   and (\ref{KP2})  in the  form
\begin{align}
\label{WAt}
&W_{At}+\left(\dfrac{V}{3}+\alpha q^2\right)W_{Ax}-\dfrac{1}{3}W_A V_x+2\alpha qDW_{A}=0,\\
\label{WBt}
&W_{Bt}+\left(\dfrac{V}{3}+\alpha q^2\right)W_{Bx}+W_A\dfrac{B_x}{6}+\alpha (W_B \,Dq+2q\,DW_{B}+W_A\, Dp)=0,\\
\label{WVt}
&W_{Vt}+\left(\dfrac{V}{3}+\alpha q^2\right)W_{Vx}-\dfrac{1}{3}W_{A}A_x+2\alpha(W_V\,Dq+q\, DW_{V}+W_B\, Dp)=0,\\
\label{qt}
&q_t+\left(\dfrac{V}{3}+\alpha q^2\right)q_x+\dfrac{1}{3}(V_y-qV_x)+2\alpha qDq=0,\\
\label{pt}
&p_t+\left(\dfrac{V}{3}+\alpha q^2\right)p_x+Dr=0,\\
\label{gamma}
&\left(\dfrac{B}{6}-r\right)_x+\alpha(\dfrac{W_B}{W_A}Dq+p_y+qp_x)=0,
\end{align}
where 
\[
 D:=\dfrac{\partial}{\partial y}-q\dfrac{\partial }{\partial x}.\]
 The constraints   (\ref{KP3}) can  be  written, after using (\ref{q}) in the form 
\begin{equation}
\label{const_final0}
q_x=\dfrac{k_y}{k}-q\dfrac{k_x}{k},\;\;\;p_x=c_{1y}-(qc_1)_x.\;\;\;
\end{equation}
Equations (\ref{WAt})-(\ref{WVt}),    with $p=0$ and the consistency  conditions   \eqref{KP1}-\eqref{KP3}  have been obtain \cite{Infeld92}.

We observe that equations (\ref{WAt}), (\ref{WBt}) and  (\ref{WVt}),  for 
$\alpha=0$ are identical to the Whitham modulation equations  for the KdV equation \cite{Whitham}.  Furthermore, for $\alpha=0$  equation (\ref{gamma}) can be solved exactly giving $r=B/6$.   If we assume that $A$, $B$ and $V$ are $y$-independent, we get the further integrals $q=6 h(y) W_A$ and $p=6h(y) W_B$ for a  function $h(y)$.
   Whitham was able to reduce  (\ref{WAt}), (\ref{WBt}) and  (\ref{WVt})   for $\alpha=0$     to diagonal form.
Using $e_1$, $e_2$ and $e_3$  defined in (\ref{ei}) as independent variables,  equations 
 (\ref{WAt}), (\ref{WBt}) and  (\ref{WVt})  for $\alpha=0$  take the form
\begin{equation}
\label{whitham00}
\dfrac{\partial}{\partial t}e_i+\sum_{k=1}^3\sigma_i^k\dfrac{\partial}{\partial x}e_k=0,\quad i=1,2,3,
\end{equation}
where the  matrix $\sigma_i^k$ given by 
\begin{equation}
\label{sigma}
\sigma=\dfrac{1}{3}VI-\dfrac{W_A}{6}\begin{pmatrix}
\partial_{e_1}W_A&\partial_{e_2}W_A&\partial_{e_3}W_A\\
\partial_{e_1}W_B&\partial_{e_2}W_B&\partial_{e_3}W_B\\
\partial_{e_1}W_V&\partial_{e_2}W_B&\partial_{e_3}W_V
\end{pmatrix}^{-1}
\begin{pmatrix}
2&2&2\\
e_2+e_3&e_1+e_3&e_1+e_2\\
2e_2e_3&2e_1e_3&2e_1e_2\end{pmatrix},
\end{equation}
where $I$ is the identity matrix and  $\partial_{e_i}W_A$ is the partial derivative with respect to $e_i$ and the same notation holds for the other quantities.
Equations (\ref{whitham00}) is a system of quasi-linear equations for $e_i=e_i(x,t)$, $j=1,2,3$. Generically, a quasi-linear  $3\times 3$ system cannot be reduced to a  diagonal form. However  Whitham, analyzing the form of the matrix $\sigma$,   was able to get the Riemann invariants that reduce the system to diagonal form. Indeed by making   the change of coordinates  
\begin{equation}
\label{Riemann}
\beta_1=\dfrac{e_2+e_1}{2},\;\;\beta_2=\dfrac{e_1+e_3}{2},\;\;\beta_3=\dfrac{e_2+e_3}{2},
\end{equation}
with  $\beta_3<\beta_2<\beta_1,$
and introducing  a matrix $\mathcal{C}$ that produces the change of coordinates $(\beta_1,\beta_2,\beta_3)^t={\cal C} (e_1,e_2,e_3)^t$,
the velocity matrix $\sigma $   in(\ref{sigma}) transforms  to  diagonal form 
\[
\tilde{\sigma}=\mathcal{C} \sigma\mathcal{ C}^{-1}=\begin{pmatrix}
v_1&0&0\\0&v_2&0\\0&0&v_3
\end{pmatrix},
\]
where the speeds $v_i=v_i(\beta_1,\beta_2,\beta_3)$ have been calculated by Whitham \cite{Whitham} and take the form (\ref{vi}).
Summarizing, the Whitham modulation equations  for KdV in the dependent variables $\beta_1>\beta_2>\beta_3$ take  the diagonal form
\[
\dfrac{\partial }{\partial t}\beta_i+v_i(\beta_1,\beta_2,\beta_3)\dfrac{\partial }{\partial x}\beta_i=0,\quad i=1,2,3.
\]

Using the same change of variables for the  first three equations  (\ref{WAt})-(\ref{WVt})  in the Whitham system for KP, this gives after  similar  computations (done  in a straitforward way with Maple)
the system of equations 
\begin{align}
\label{eqb}
&\dfrac{\partial}{\partial t}\beta_i+(v_i+\alpha q^2)\dfrac{\partial }{\partial x}\beta_i+\alpha(2 q  D\beta_i-(v_i-2\beta_i)Dq+Dp)=0,\quad i=1,2,3,\\
\label{eqq1}
&\dfrac{\partial}{\partial t}q+\left(
\dfrac{1}{3}\sum_{i=1}^3\beta_i+\alpha q^2\right)q_x+2\alpha 
Dq+\dfrac{1}{3}D(\sum_{i=1}^3\beta_i)=0,\\
\label{eep}
&p_t+\left(
\dfrac{1}{3}\sum_{i=1}^3\beta_i+\alpha q^2\right)p_x+Dr=0,\\
\label{eer}
&r_x-\dfrac{B_x}{6}-\alpha(c_1(q_y-qq_x)+p_y+qp_x))=0,
\end{align}
with $D=\dfrac{\partial}{\partial y}-q\dfrac{\partial}{\partial x}$ 
and the  quantities $B$ and $c_1$ take the form
\[
B=\sum_{i=1}^3\beta_i^2-2(\beta_1\beta_2+\beta_2\beta_3+\beta_1\beta_3), \quad c_1=\beta_3+\beta_2-\beta_1+2(\beta_1-\beta_3)\dfrac{E(s)}{K(s)}
\]

The constraints (\ref{const_final0}) can be written in the dependent  variables $\beta_1>\beta_2>\beta_3$,  in the form
\begin{equation}
\label{const_final}
q_x=\sum_{j=1}^3\dfrac{\beta_{jy}-q\beta_{jx}}{3v_j-V},\quad p_x=\sum_{j=1}^3\dfrac{2\beta_j-V}{3v_j-V}\left[\beta_{jy}-q\beta_{jx}\right].\;\;\;
\end{equation}

The  equations  (\ref{eqb}) (\ref{eqq1}) and (\ref{const_final}) are equivalent to the equations obtained  in \cite{Biondini}, while the equations (\ref{eep}) and (\ref{eer}) are new.  For example the  equation  for the variable $p$  in \cite{Biondini}
 is  the  linear combination of the two constraints (\ref{const_final}), namely
 \[
 p_x+(\beta_1+\beta_3-\beta_2)q_x=\dfrac{E(s)}{K(s)}D\beta_1+\left(1-\dfrac{E(s)}{K(s)}\right)D\beta_3.
 \]
 \begin{remark}
For $\alpha=0$  equation (\ref{eer}) can be solved exactly giving   the integral $r=B/6+g$, where $g=g(y,t)$ is an arbitrary function.
If we further assume that $\beta_i$, $i=1,2,3$  are $y$-independent,  and  we set $g(y,t)=0$, then  we get the  integrals $q=6 h(y) W_A$ and $p=-6h(y) W_B$ for an arbitrary  function $h(y)$ and the equations (\ref{eqb})   coincide with the Whitham modulation equations for the KdV equation.
If we assume like in \cite{Biondini} that the quantities $\beta_i(x,y,t)=\beta_i(\eta,t)$, $i=1,2,3$,  where $\eta=x+P(y,t)$ and $q=P_y(y,t)$, $r=r(\eta,t)$ and $p=p(y,t)$ one obtains
$D\beta_i=0$ and $Dr=0$ and the Whitham-KP system reduce to 
\begin{align}
\label{eqb_r}
&\dfrac{\partial}{\partial t}\beta_i+(v_i+\alpha q^2)\dfrac{\partial }{\partial \eta}\beta_i-\alpha((v_i-2\beta_i)q_y+p_y)=0,\quad i=1,2,3,\\
\label{eqq1_r}
&\dfrac{\partial}{\partial t}q+2\alpha 
qq_y=0,\quad p_t=0,\quad r_\eta-\dfrac{B_\eta}{6}-\alpha(c_1q_y+p_y)=0.
\end{align}
 In equation (\ref{eqb_r})  and the second equation in (\ref{eqq1_r}),   since $\beta_i=\beta_i(\eta,t)$  and $r=r(\eta,t)$,   namely they are independent from $y$,  consistency conditions imply that  $q_y=0$ or $q_y=const$ and $p_y=0$ or $p_y=const$ which give the reduction to KdV or cylindrical KdV. For further details refer to \cite{Biondini}. 
\end{remark}
%

\subsection{ Limiting behaviour of the  Whitham modulation equations  near the soliton front}
In  the limit $m\to 1$ the wave-train of oscillations  becomes a sequence of near-solitary waves.
When $m\to 1$  one has  (see e.g. \cite{Lawden})
\begin{equation}
\label{ellipticlead}
E(m)\simeq 1+(1-m)\left[\Lambda-\frac{1}{2}\right],
\quad K(m)\simeq \Lambda,\quad \Lambda=\dfrac{1}{2}\log\dfrac{16}{1-m^2}.
\end{equation}
 One can verify that the speeds $v_i$ have the following limiting 
 behaviour   ( see e.g.  \cite{GP}) in the `solitonic limit'',  $m=1$ or  $\beta_2=\beta_1$:
\begin{equation}
\begin{split}
\label{limit_edge}
&v_1(\beta_1,\beta_1,\beta_3)=v_2(\beta_1,\beta_1,\beta_3)=\dfrac{2}{3}\beta_1+\dfrac{1}{3}\beta_3,\\
&v_3(\beta_1,\beta_1,\beta_3)=\beta_3.
\end{split}
\end{equation}
In this limit  the equation for  the variable  $\beta_3$  in (\ref{eqb})  takes the form 
\[
\dfrac{\partial}{\partial t}\beta_3+\beta_3\dfrac{\partial }{\partial x}\beta_3+\alpha\left( (q\beta_3+p)_y+q( \beta_{3y}-(q\beta_3+p)_x)\right)=0.
\]
This equation has to be equivalent to the dKP equation (\ref{dKP}). Indeed using the  linear combination of the  constraints (\ref{const_final}) one obtains, in the limit $\beta_2\to \beta_1$, the equation
$
p_x+\beta_3 q_x=D\beta_3$
which  implies  the dKP equation
\[
\dfrac{\partial}{\partial t}\beta_3+\beta_3\dfrac{\partial 
}{\partial x}\beta_3+\alpha (q\beta_3+p)_y=0,\quad 
 \beta_{3y}-(q\beta_3+p)_x=0,
\]
or
\[
\left(\dfrac{\partial}{\partial t}\beta_3+\beta_3\dfrac{\partial 
}{\partial x}\beta_3\right)_x+\alpha\dfrac{\partial^2}{\partial y^2}\beta_3=0.
\]
The above equation implies that if we chose  $\beta_3(x,y, 0)=0$ at the soliton front,   it will remain zero also al later times.
It follows that  when $m\to 1$ and $\beta_3=0$, we have  $p_x=0$, $p_y=0$, $B=0$, and $c_1=0$ so that the Whitham system reduce to the form
\begin{align}
\label{soli1}
& \beta_{1t}+\left(\frac{2}{3}\beta_1-\alpha q^2\right)\beta_{1x}+2\alpha q\beta_{1y}+\frac{4}{3}\beta_1\alpha (q_y-qq_x)=0,\\
 \label{soli2}
 &q_t+(\frac{\beta_1}{3}-\alpha q^2)q_x+2\alpha qq_y+\frac{2}{3}(\beta_{1y}-q\beta_{1x})=0\\
&p_t+r_y=0,\quad r_x=0,
\end{align}
namely we have two sets of uncouple equations.
 It is straightforward to check that the first two equations of the above system are elliptic (see below). In the next section we want to show that the  equations
  (\ref{soli1}) and (\ref{soli2}) can be derived as modulation of the soliton parameters.

\subsection{Soliton modulation of the KP equation}
We are interested in studying the slow modulation of the wave 
parameters of the soliton solution (\ref{soliton}) following 
Whitham's averaging procedure of the Lagrangian density.
We make the ansatz
\[
\psi_x=a\,\mbox{sech}^2\left[\left(\dfrac{a}{12}\right)^{\frac{1}{2}}\left(x-\frac{\omega}{k} t +\frac{l}{k}y\right)\right],\quad \psi_t=-\dfrac{\omega}{k}\psi_x,\;\;\;\;\psi_y=\dfrac{l}{k}\psi_x,
\]
where  $a$ is the amplitude, $k$ the wave number and $\omega$ the frequency.
 The average Lagrangian $\mathcal{L}$  is obtained by integration, namely
  \begin{equation}
 \mathcal{L}=k\int_{-\infty}^{+\infty}L 
 \d x=\dfrac{4}{15}\sqrt{12}\left(ka^{\frac{5}{2}}-5a^{\frac{3}{2}}\omega+5\alpha\frac{l^2}{k}a^{\frac{3}{2}}\right).
 \end{equation}
%
The variation with respect to the amplitude gives
\begin{equation}
\label{omega_soliton}
\dfrac{\delta {\cal L}}{\delta a}=0\quad \longrightarrow  \quad \omega=\dfrac{ka}{3}+\alpha \dfrac{l^2}{k}.
 \end{equation}
 The variation with respect to the phase $\theta(x,y,t)=kx+ly-\omega t$ gives the equations
 \[
 \dfrac{\partial}{\partial x}\dfrac{\delta {\mathcal L}}{\delta  k}-\dfrac{\partial}{\partial t}\dfrac{\delta {\mathcal L}}{\delta \omega}+\dfrac{\partial}{\partial y}\dfrac{\delta {\mathcal L}}{\delta l}=0,
 \]
 namely
 \begin{equation}
 \label{eq_a}
 a_t+\left(\frac{a}{3}-\alpha q^2\right)a_x+\frac{4}{3}a\alpha 
 (q_y-qq_x)+2\alpha qa_y=0,
 \end{equation}
plus  the consistency equations
 \[
\dfrac{\partial}{\partial y}k-\dfrac{\partial}{\partial x}{l}=0,\quad \dfrac{\partial}{\partial t}k+\dfrac{\partial}{\partial x}{\omega}=0,\quad   \dfrac{\partial}{\partial t}l+\dfrac{\partial}{\partial y}{\omega}=0, \]
that can be written in the form
 \begin{align}
&k_y=(qk)_x,\quad q=\dfrac{l}{k},\\
\label{eq_k}
& k_t+(\frac{a}{3}-\alpha q^2)k_x+2\alpha qk_y+\frac{k}{3}a_x=0,\\
\label{eq_q}
&q_t+(\frac{a}{3}-\alpha q^2)q_x+2\alpha qq_y+\frac{1}{3}(a_y-qa_x)=0.
 \end{align}
We have three equations (\ref{eq_a}), (\ref{eq_k}) and (\ref{eq_q})  for three variables $a,k$ and $q$, while $\omega$ is recovered from (\ref{omega_soliton}).
The   equations (\ref{eq_a})  and (\ref{eq_q}) are independent from the variable $k$, 
\[
\begin{pmatrix}
a\\q\end{pmatrix}_t+\begin{pmatrix}\frac{a}{3}-\alpha q^2&-\frac{4}{3}\alpha qa\\
-\frac{q}{3}&\frac{a}{3}-\alpha q^2\end{pmatrix}\begin{pmatrix}
a\\q\end{pmatrix}_x+
\begin{pmatrix}
2\alpha q&\frac{4}{3}\alpha a\\
\frac{1}{3}&2\alpha q\\
\end{pmatrix}\begin{pmatrix}
a\\q\end{pmatrix}_y=0.
\]
Defining  $A_1$  as the first matrix and $A_2$ as the second matrix, the  above  system of equations is strictly hyperbolic if 
the eigenvalues of 
\[
A_1+\xi A_2
\]
are real for any real $\xi$.
 After a simple calculation one obtains that the eigenvalues $\lambda_i$, $i=1,2$, of the matrix $A_1+c\xi A_2$ are 
 \[
\lambda_{1,2}=\dfrac{a}{3}-\alpha q^2+2\xi\alpha q\pm \dfrac{2}{3}\sqrt{\alpha a(q-\xi)^2},
\]
where the  amplitude  $a>0$. From the above expression,  it is clear that for KPII  ($\alpha=1$) all the eigenvalues are always real while
 for KPI  ($\alpha=-1$) the eigenvalues 
are complex. In this case it is expected that the parameters describing the evolution of the leading soliton  front  have a singularity of elliptic type  (elliptic umbilic catastrophe) as in the 
singularity formation of the semiclassical limit of the nonlinear 
Schr\"odinger equation.  Indeed  in this  case
the generic initial data evolve, near the point of elliptic umbilic catastrophe,   into a breather, that is a rational solution.
For the KPI case, we numerically observe that the  leading solitons emerging from the dispersive shock wave always break into a series 
of lumps arranged on a lattice.

\section{Solutions to  focusing  NLS  and KPI   equations}
The Cauchy problem for the semiclassical limit of the
focusing NLS equation
\begin{equation}
\label{NLS}
\ii\epsilon\psi_y+\dfrac{\epsilon^2}{2}\psi_{xx}+\psi|\psi|^2=0,
\end{equation}
where we denote time by $y$,
was considered in  \cite{KMM}. For generic initial data  
$\psi(x,y=0;\epsilon)$ the  solution develops an oscillatory zone. 
The  $(x,y)$ plane is basically divided into two regions,  a region where the solution $\psi(x,y;\epsilon)$ has a highly oscillatory behaviour  with oscillations of wave-length $\epsilon$, and a region where the solution
is non oscillatory.  In  \cite{DGK} and \cite{Bertola_Tovbis},  the transition region between these two regimes has been considered.
Introducing the slow variables 
\[
\rho=|\psi|^2,\quad w=\dfrac{\epsilon}{2 \ii}\left(\dfrac{\psi_x}{\psi}-\dfrac{\overline{\psi}_x}{\overline{\psi}}\right),
\]
the NLS equation can be written in the form
\begin{align}
&\rho_y+(\rho w)_x=0,\\
&w_y-\rho_x+ww_x+\dfrac{\epsilon^2}{4}\left(\dfrac{\rho_x^2}{2\rho^2}-\dfrac{\rho_{xx}}{\rho}\right)_x=0.
\end{align}
The semiclassical limit takes the hydrodynamic form
\begin{align}
\label{eq0}
&\rho_y+(\rho w)_x=0,\\
\label{eq01}
&w_y-\rho_x+ww_x=0.
\end{align}
For generic initial data, the solution of the above elliptic system 
of equations  develops a point  $(x_0,y_0)$ where the gradients 
$\rho_x$ and $w_x$ are divergent but the quantities $w(x_0,y_0)$ and $\rho(x_0,y_0)$  remain finite.
Such a point is called an {\it elliptic umbilic catastrophe}.
Correspondingly the solution of the NLS equation remains smooth and can be approximated by the tritronqu\'ee solution to the Painlev\'e I equation  $f_{zz}=6f^2-z$ \cite{DGK}. 
However the approximation is not valid near the poles  $z_p$ of the 
tritronqu\'ee solution. At the poles   the NLS solution is approximated  \cite{Bertola_Tovbis} by the rational Peregrine breathers.
These breathers are parametrized by  the two  real  constants $a$ and $b$
 and take the form
 \begin{equation}
 \label{breather}
 Q(x,y;a,b)=\ee^{-\ii(ax+(a^2/2-b^2)y)}b\left( 1-4\dfrac{1+2\ii b^2y}{1+4b^2(x+a y)^2+4b^4 y^2/4}\right),
 \end{equation}
 where $|Q(x,y;a,b)|\to b$ as $|x|\to\infty$ and  the maximum value of $|Q(x,y;a,b)|$ is three times the background value $b$, namely
 \[
 \sup_{x\in\mathbb{R}, y\in \mathbb{R}^+}|Q(x,y;a,b)|=3b.
 \]
 
 Identifying $a=-w(x_0,y_0)$ and $b=\sqrt{\rho(x_0,y_0)}$, the NLS solution is given in the limit $\epsilon\to 0$ by \cite{Bertola_Tovbis}
 \[
 \psi(x,y;\epsilon)=\ee^{\frac{\ii}{\epsilon}\Phi(x_p,y_p)}Q\left(\dfrac{x-x_p}{\epsilon},\dfrac{y-y_p}{\epsilon}\right)+O(\epsilon^{\frac{1}{5}})
 \]
 where $\Phi(x_p,y_p)$ is a phase,
 $(x_p,y_p)$  is related to the poles $z_p$ of the tritronqu\'ee solution  $f(z)$  via the variable
 \begin{equation}
 \label{zp}
 z_p=\dfrac{c_0}{\epsilon^{\frac{4}{5}}}[x_p-x_0+(a+\ii b)(y_p-y_0)],
 \end{equation}
  \begin{wrapfigure}{r}{0.5\textwidth}
    \includegraphics[width=0.45\textwidth]{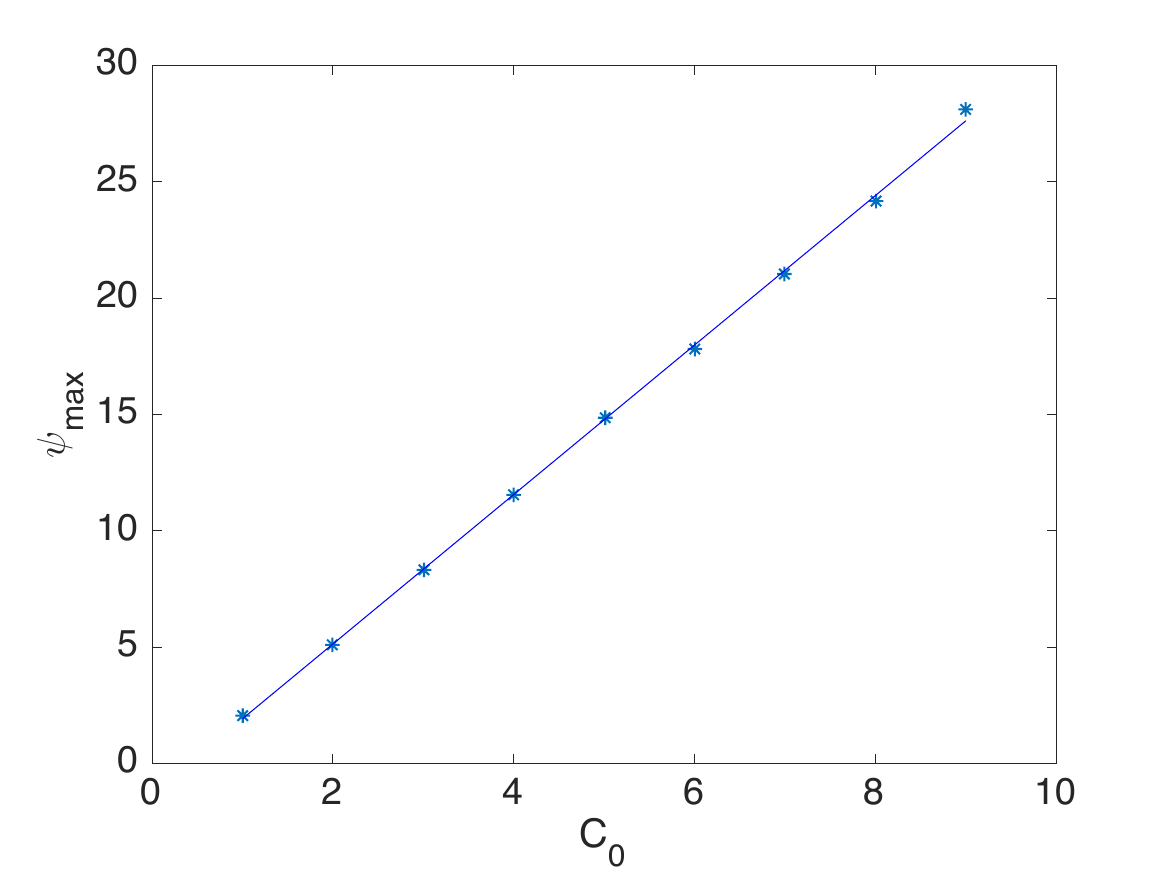}
    \label{NLSfig}
 \caption{ The maxima $\psi_{max}$ of the $L^{\infty}$ norm of the solution to the focusing NLS equation  (\ref{NLS}) 
 for the   initial data $C_{0}\mbox{sech}^{2}x$ for 
  several values of $C_{0}=1,\ldots,9$  with a linear fit  $|\psi(x,y;\epsilon)|_{max}=3.2128 C_{0}
   -1.2864$.} 
 \end{wrapfigure}

 with $(x_0,y_0)$ the point of elliptic umbilic catastrophe and $c_0$ 
 a constant that depends on the initial datum.  For example the first
 breather corresponds to the first  pole  at $z_p\simeq-2.38$  on the negative real axis of the 
 tritronqu\'ee solution.   
The macroscopic feature of this behaviour is that  the maximum hight 
of the solution is approximately $3$ times the value $b$   that is 
the value of  $\rho(x_0,t_0;\epsilon=0)$ at the critical point. 
Furthermore the  above formula for $z_p$ shows that the position of the lump in the $(x,t)$  plane   scales like $\epsilon^{\frac{4}{5}}$.
When  the value $b$ is not available,  one may wonder whether the maximum peak of the NLS solution scales linearly with 
the maximum value of the initial data. Using the same numerical approach as 
in \cite{DGK} (we use $N=2^{14}$ Fourier modes and $N_{t}=10^{4}$ 
time steps), we get  the $L^{\infty}$ norm  of $\psi(x,y;\epsilon)$ for $\epsilon=0.1$ and  for the initial data 
$\psi(x,0)=C_{0}\partial_{x}\mbox{sech}^{2}x$ for several values of   
$C_{0}$.
The maxima of the 
$L^{\infty}$ norms are shown in Fig.~\ref{NLSfig} 
in dependence of $C_{0}$. They can be fitted via linear regression to the line
$3.2128 C_{0}
   -1.2864$, thus confirming that the maximum value of the solution 
   scales linearly with  the maximum value of the initial data above some threshold amplitude $C_0$.

 We  now connect the NLS breather solution (\ref{breather})  to the KPI lump solution  (\ref{lumpgeneral})  by observing    that the expression
\begin{equation}
\label{lump_breather}
u(x,y,t)=12\left|Q\left(x-(a^2+3b^2)t,2\sqrt{3}(y+2at); \frac{a}{2\sqrt{3}},\frac{b}{2}\right)\right|^2-3 b^2
\end{equation}
coincides with the general lump solution  (\ref{lumpgeneral}) of KPI   for $\epsilon=1$.
Using this connection, we make the following conjecture.
\begin{conjecture}
The position of the lumps emerging from the soliton front is   determined by the relation
 \begin{equation}
 \label{zp_lump}
 z_p=\dfrac{c_0}{\epsilon^{\frac{4}{5}}}[x_p-(a^2+3b^2)t_p-(x_0-(a^2+3b^2)t_0)+(a+\ii\sqrt{3}b)(y_p+2at-y_0-2at_0)],
 \end{equation}
where $(x_0,y_0)$ is the position where a singularity of the Whitham system is expected to appear at the time $t_0$ and 
$(x_p,y_p)$ is the position where the lump is expected to appear at the time $t_p$ and $c_0$, $a$ and $b$ are some constants.
\end{conjecture}
For initial data symmetric with respect to  $y\to -y$   the first lump that is appearing in the KPI solution is on the line $y=0$, thus  
 $a=0$ and $y_0=y_p=0$ due to symmetry reasons.
 We conclude form (\ref{zp_lump}) that the position of the first lump is expected to be given by 
  \begin{equation}
 \label{zp_lumps}
 z_p=\dfrac{c_0}{\epsilon^{\frac{4}{5}}}[x_p-3b^2t_p-(x_0-3b^2t_0)],
 \end{equation}
 namely the quantity $x_p-3b^2t_p$  is expected to  scale like $\epsilon^{\frac{4}{5}}$. We are going to verify this ansatz numerically in the next section.

%
%
\section{Numerical Method}
\label{sec:numerical_method}
In this section we summarize the numerical methods used in the 
following section to solve the Cauchy problem for KPI in the small 
dispersion limit. 
We consider the evolutionary form of the KP equation (\ref{KP}):
\begin{equation}
  u_t+uu_x+\epsilon^2 u_{xxx}=-\alpha\partial_x^{-1}u_{yy},
  \label{eq:KPI_num}
\end{equation}
defined on the periodic square $[-5\pi,5\pi]^2$, with initial 
condition $u(x,y,0)=u_{0}(x,y)$; here $\partial_{x}^{-1}$ is defined 
via its Fourier multiplier $-\mathrm{i}/k_{x}$ where $k_{x}$ is the dual 
Fourier variable to $x$.


For the numerical approximation of the solution $u(x,y,t)$ of equation~\eqref{eq:KPI_num}, we adopt a \emph{Fourier collocation} method (also known as \emph{Fourier pseudospectral} method) in space coupled with a \emph{Composite Runge--Kutta} method in time.

Referring to~\cite{CHQZ2,trefethen} for a detailed overview of Fourier collocation methods and spectral methods in general, we sketch here the main features of this discretization method.
The starting point of Fourier spectral methods consists in 
approximating the Fourier transform $\widehat{u}(k_{x},k_{y},t)$ of the 
solution $u(x,y,t)$, where $k_{x}$, $k_{y}$ are the dual variables to 
$x$, $y$, via a discrete Fourier transform for which fast 
algorithms exist, the \emph{fast Fourier transform} (FFT). This means we 
approximate the rapidly decreasing initial data as a periodic (in $x$ 
and $y$) function. We will always work on the domain 
$5[-\pi,\pi]\times5[-\pi,\pi]$ in the following.  We use $N_{x}$ 
respectively $N_{y}$ collocation points in $x$ respectively $y$.

The discretized approximation of the KPI equation~\eqref{eq:KPI_num} can be written in the form:
\begin{equation}
  \widehat{u}_t=\mathbf{L}\widehat{u}+\mathbf{N}(\widehat{u}),
  \label{eq:ut_LN}
\end{equation}
where for the KPI equation~\eqref{eq:KPI_num}, the linear and nonlinear parts $\mathbf{L}$ and $\mathbf{N}$ have the form:
\begin{equation}
  \begin{aligned}
  &\mathbf{L} = -\ii\frac{k_y^2}{k_x}+\epsilon^2\ii k_x^3, \\
  &\mathbf{N}(\widehat{u}) = -\frac12\ii k_x\widehat{u^2}.
  \end{aligned}
  \label{eq:def_LN}
\end{equation}
The convolution in Fourier space in the nonlinear term $\mathbf{N}$ 
in equation~\eqref{eq:def_LN} is computed in physical space followed 
by a two-dimensional FFT.

For the time discretization of equation~\eqref{eq:ut_LN} 
several fourth order methods were discussed in \cite{KR11} for the 
small dispersion limit of KP. We adopt here Driscoll's Composite 
Runge--Kutta method~\cite{Dri}, which requires that the linear 
operator $\mathbf{L}$ of equation~\eqref{eq:def_LN} is diagonal, which 
is the case here. Thus the evaluation of both positive and negative powers of $\mathbf{L}$ can be obtained with a computational cost $O(N)$.

Composite Runge--Kutta methods partition the Fourier space for the 
linear part of the equation into two parts, one for the low frequencies (or ``slow'' modes), $|\mathbf{k}|:=|(k_x,k_y)|<k_{\mathrm{cutoff}}$, and one for the high frequencies (or ``stiff'' modes), $|\mathbf{k}|\ge k_{\mathrm{cutoff}}$.
Then, the Fourier components of the solution are advanced in time 
using different Runge--Kutta integrators for the two partitions. In 
particular, a third-order $L$-stable method (RK3 in the following) is 
used for the higher frequencies, while for the lower frequencies 
stiffness is not an issue and a standard explicit fourth-order method 
(RK4) can be used. As a result, the method is explicit, but has much 
better stability properties than the explicit RK4 method for which no 
convergence could be observed in the studied examples in \cite{KR11}. 
Despite the use of a third order method for the high frequencies, 
Driscoll's method shows in practice fourth order accuracy as shown in 
\cite{KR11} and references therein.

In his article~\cite{Dri}, Driscoll suggests to adopt the fourth order method for all the frequencies such that:
\begin{equation}
  ||\mathbf{L}||_{\infty} < \frac{2.8}{h},
  \label{eq:driscoll_cutoff}
\end{equation}
where $h$ is the time-step used, in accordance with the stability 
region of the RK4 method, see e.g.~\cite{trefethen}.
However, in previous studies as \cite{KR11} and references therein, 
it was observed that the method is stable only  if very  small time 
steps  depending on the spatial resolution are used (obviously 
$h\propto 1/(N_{x}N_{y})$).
For this reason, we modified condition~\eqref{eq:driscoll_cutoff} to the following:
\begin{equation}
   ||\mathbf{L}||_{\infty} < \frac{2^{-7}}{h}.
  \label{eq:my_cutoff}
\end{equation}
As a result of this change, many fewer Fourier modes of the 
linear part are 
advanced with a RK4 method  than in Driscoll's original method, but 
  this is still preferable over a standard RK3 method (an 
explicit RK3 method would impose similar stability requirement as 
RK4, and an implicit method would make the solution of an implicit 
equation system necessary in each time step, which would be 
computationally too expensive).

Due to the very high accuracy required by our simulations, the numerical method exposed so far has been implemented in a MPI-parallel \texttt{C} code.

The accuracy of the solutions is controlled as in \cite{KR11} in two 
ways: since the KPI solution for smooth initial data is known to stay 
smooth, its Fourier transformed must be rapidly decreasing for all 
time. Thus if the computational domain is chosen large enough, this must 
be also the case for the discrete Fourier transform. The decrease of 
the Fourier coefficients can thus be used to control the numerical 
resolution in space during the computation. If the latter is assured, 
the resolution in time can be controlled via conserved quantities of 
the KP solution as the $L^{2}$ norm or the energy, computed numerically as:
\begin{equation}
  \|u\|_{L^2}^2 = \sum_{|\mathbf{k}|=0}^N|\widehat{u}_{\mathbf{k}}|^2,
  \label{eq:L_2_plancherel}
\end{equation}
which will be 
numerically time dependent due to unavoidable numerical errors. As 
discussed for instance in \cite{KR11} the accuracy in the 
conservation of such quantities can be used as an indicator of the 
numerical accuracy.



\section{Numerical solution}
In this section we analyse the behaviour of the KPI solution for the initial data
\begin{equation}
\label{initial}
 u_0(x,y)=-C_0\partial_x\sech^2\sqrt{x^2+y^2}.
 \end{equation}
for several values of $\epsilon$ and $C_0$.
In table~\ref{tab:parameters}, we report the different set-ups for the numerical simulations.
\begin{table}
  \caption{Parameter values for the numerical experiments run 
  (numbered by $n$) in this work.}
  \label{tab:parameters}
  \centering
  \begin{tabular}{lcccc}
    \toprule
    n & $\epsilon$ & $C_0$ & $h$ & grid \\
    \midrule
    1 & $0.02$     & 6 & $4\cdot10^{-5}$ & $2^{15}\times2^{15}$ \\
    2 & $0.03$     & 6 & $4\cdot10^{-5}$ & $2^{15}\times2^{15}$ \\
    3 & $0.04$     & 6 & $8\cdot10^{-5}$ & $2^{15}\times2^{15}$ \\
    4 & $0.05$     & 6 & $1\cdot10^{-4}$ & $2^{15}\times2^{15}$ \\
    5 & $0.06$     & 6 & $1\cdot10^{-4}$ & $2^{14}\times2^{14}$ \\
    6 & $0.07$     & 6 & $2\cdot10^{-4}$ & $2^{14}\times2^{14}$ \\
    7 & $0.08$     & 6 & $1\cdot10^{-4}$ & $2^{14}\times2^{14}$ \\
    8 & $0.09$     & 6 & $2\cdot10^{-4}$ & $2^{14}\times2^{14}$ \\
    9 & $0.10$     & 6 & $2\cdot10^{-4}$ & $2^{14}\times2^{14}$ \\
    10& $0.10$     & 4 & $1\cdot10^{-4}$ & $2^{13}\times2^{13}$ \\
    11& $0.10$     & 5 & $1\cdot10^{-4}$ & $2^{13}\times2^{13}$ \\
    12& $0.10$     & 7 & $1\cdot10^{-4}$ & $2^{13}\times2^{13}$ \\
    13& $0.10$     & 8 & $1\cdot10^{-4}$ & $2^{13}\times2^{13}$ \\
    \bottomrule
  \end{tabular}
\end{table}

The solution $u(x,y,t;\epsilon)$ starts to oscillate  around the time  and the location where  the solution of the 
dKP (\ref{dKP}) equation has its first  singularity, which occurs on the positive part of the initial data.
There is a second singularity that occurs slightly later on the negative part of the initial data and a  dispersive shock wave develops also  there.
The  two dispersive shock wave fronts behave quite differently in time.
While in the negative front the oscillations are defocused, 
in the positive front the oscillations seem to be focused 
and the (modulated) line soliton fronts  break  into a number of   lumps that are arranged in 
a  lattice, as shown in Fig.~\ref{soliton_lattice}.
 \begin{figure}[htbp]
\centering
   \includegraphics[width=0.9\textwidth]{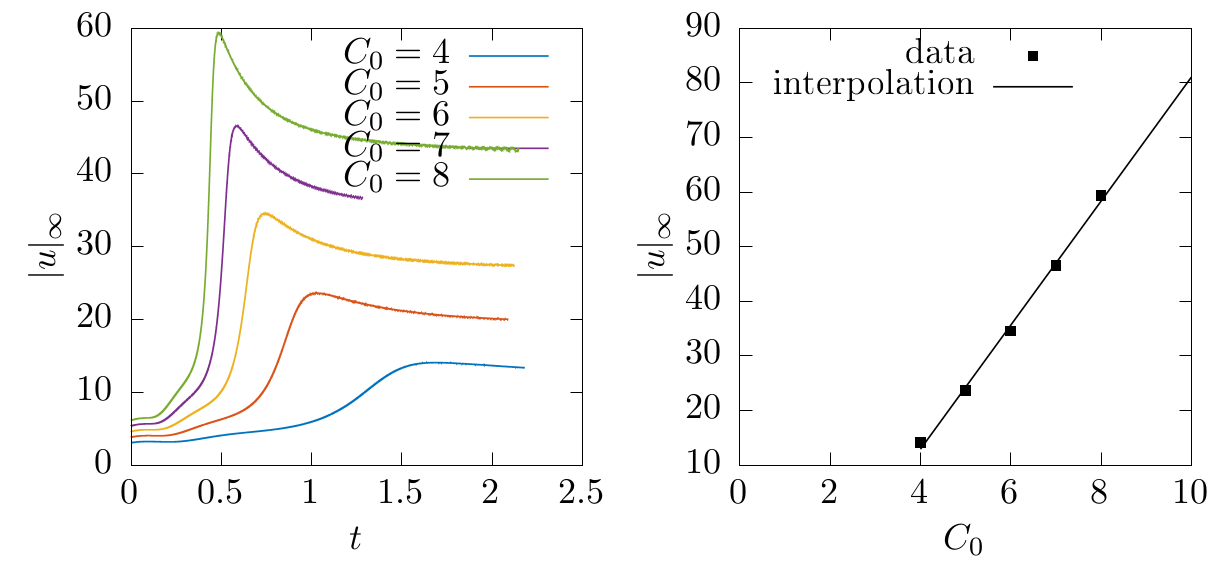}
 \caption{The $L^{\infty}$  norm  the solution of the  KPI equation  as a function of time for $\epsilon=0.10$   and for different values of the initial amplitude. The interpolation expression is $|u|_\infty=10.838C_0-29.894$.}
 \label{C0}
\end{figure}

According to a result of \cite{FS97}, for small norm initial data
\[
\int_{-\infty}^{\infty}\int_{-\infty}^{+\infty}\de y \de\xi \widehat{u}_0(\xi,y)\ll 1,
\]
 the solution of  the KPI equations with $\epsilon=1$ does not develop lumps.
 Here $\widehat{u}_0(\xi,y)$ is the Fourier transform with respect to $x$ of the initial data.
When we introduce the small $\epsilon$ parameter, such norm is of order $1/\epsilon^2$ and therefore it is never small.
For this reason, the evolution of our initial data always  develops  
lumps for sufficiently small $\epsilon$.  However for small values of $C_0$, namely for $C_0\leq 3$  the
the times required by the solution to develop the first  lump were regarded as too long,  and thus disregarded.

The same question as for NLS in Fig.~\ref{NLSfig} 
is addressed in  Fig.~\ref{C0} for the KPI 
example. We show for several values of the constant $C_0$  and for 
fixed $\epsilon=0.1$  the maximum amplitude as a function of time. The amplitude of the first lump  is proportional to   the initial amplitude.
We then consider  on the right in Fig.~\ref{C0}  the  maximum of 
the $L^{\infty}$ norm in the range of time considered   as a function of the maximum amplitude of the initial data  $u_0(x,y,0)$ in (\ref{initial})  which is proportional to  $C_0$. 
The fitting shows that $u_{max}$ is approximately  $10.8$  times $C_0$. 

\begin{figure}[tbp]
 \centering
  \includegraphics[width=\textwidth]{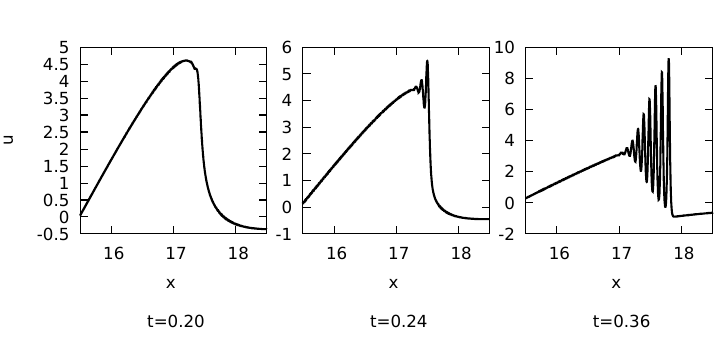}\\
 \includegraphics[width=\textwidth]{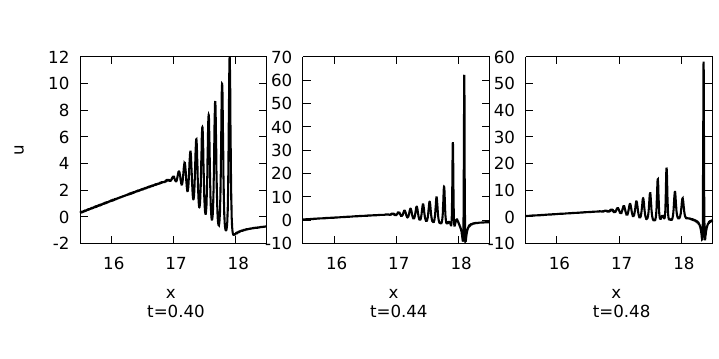}
    \caption{Solution to KPI equation for the initial data  
  $-6\partial_x\sech^2\sqrt{x^2+y^2}$ along the line $y=0$ for 
  $\epsilon=0.02$  for several values of time.  The formation of the lump and its detachment from the train of oscillations can be clearly seen.}
    \label{figure4}
 \end{figure}

\begin{figure}[htbp]
\centering
\includegraphics[width=\textwidth]{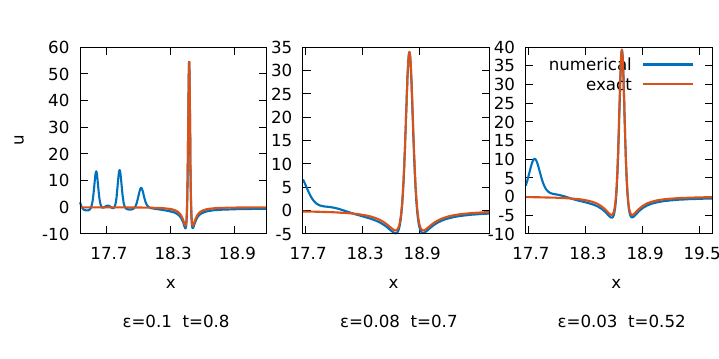}
 \caption{A comparison between the numerical solution and the lump 
 formula~\eqref{lumpgeneral} for four different values of $\epsilon$, at a time slightly after the lump achieves its maximum height.
  The maximum peak becomes   narrower  and higher with decreasing 
  values of $\epsilon$.}
  \label{figure5}
\end{figure}

In Fig.~\ref{figure4}, one can see the formation of the first lump 
from the dispersive shock of KPI on the $x$-axis.
Next we consider the fitting of the first  spike  that emerges from 
the soliton front to the KP lump (\ref{lumpgeneral}). This is shown in Fig.~\ref{figure5} 
on the $x$-axis for various values of $\epsilon$. The excellent 
agreement is obvious.

\begin{figure}[htbp]
\centering
       \includegraphics[width=\textwidth]{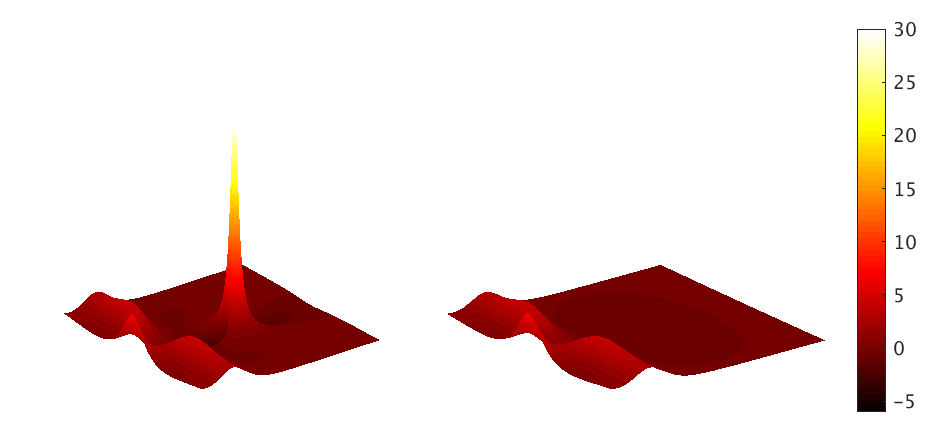}
   \caption{2D plot of the KPI solution for $\epsilon=0.06$ and $t=0.9$. On the right picture, the maximum peak has been subtracted using the lump solution (\ref{lumpgeneral}).}
   \label{figure4b}
\end{figure}

In Fig.~\ref{figure4b} we show the 2D-plot of the highest peak. 
We subtract the fitted lump solution and, as can be seen from the picture,
the difference is negligible with respect to the  remaining oscillations.

\begin{figure}[htbp]
\centering
\includegraphics[width=1.0\textwidth]{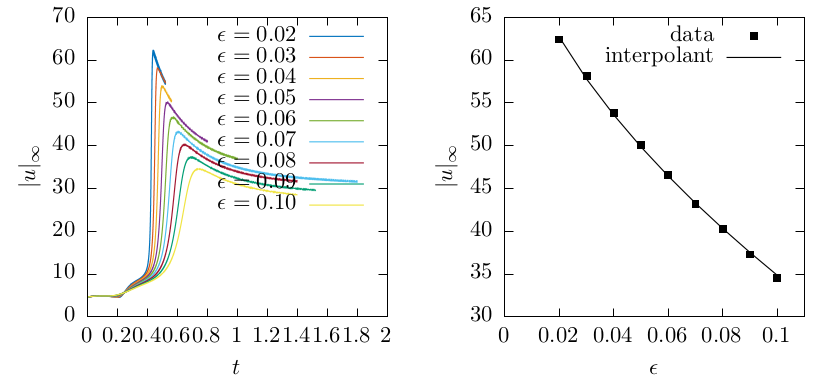}
   \caption{On the left the   $L^{\infty}$ norm  of the solution 
   $u(x,y,t;\epsilon)$  as a function of time for several   values of 
   $\epsilon$. On the right a fitting of
   $|u|_{\infty}$ in dependence of $\epsilon$ to $c_1+c_2\epsilon^\beta$, which yields $c_1 = 77.9350$, $c_2=-191.4782$, $\beta=0.6437$.
   }
  \label{figure6}
\end{figure}

We  study numerically the scaling of the lump parameters as a function of $\epsilon$ for fixed initial data.
The first scaling that we consider is the $L^{\infty}$ norm $|u|_{\infty}$ as a function of $\epsilon$  (see Fig.~\ref{figure6}).
A fitting of $|u|_{\infty}$ to $c_{1}+c_{2}\epsilon^{\beta}$ with 
gives  $c_{1}=77.9350$ $c_{2}=-191.4782 $ and   $\beta=0.6437$.

\begin{figure}[tbp]
\centering
  \input{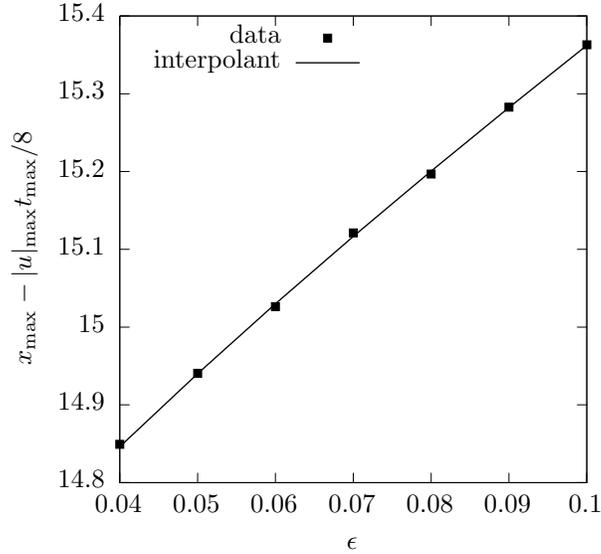}
     \caption{
     The value 
     $x_\text{max}-|u|_\text{max}/8t_\text{max}$ as a function of 
     $\epsilon$. A power fitting $x_\text{max}-|u|_\text{max}t_\text{max}/8=c_1+c_2\epsilon^\beta$ gives the coefficients $c_1=14.354$, $c_2=6.1037$, $\beta=0.7820$. 
     }
      \label{figure7}
    \end{figure}

Next we consider the dependence of the position and the time of  
appearance of the highest peak  as a function of  $\epsilon$. Since 
the time of the second breaking, its location and the value of the 
solution are not known, but all enter formula (\ref{zp_lump}), it 
will be numerically inconclusive if they all will be identified via 
some fitting for $x$, $t$  and $u$ separately. Instead we just 
consider the combination of these values needed for (\ref{zp_lump}), 
$x_\text{max}-|u|_\text{max}/8t_\text{max}$ and fit the observed 
values to $c_1+c_2\epsilon^\beta$. As shown in Fig.~\ref{figure7}, we find
$c_{1}=14.3537$,    $c_{2}=6.1037$ and     $\beta=0.7820$  which is compatible with the value $4/5$.
%

\section{Conclusion}
In this work we have presented a detailed numerical study of the long 
time behavior of dispersive shock waves in KPI solutions. It was 
shown that in the positive part of the solution, a secondary breaking 
of the dispersive shock wave can be observed for 
 sufficiently long times,  depending on the amplitude of the initial data. At this 
secondary breaking, the parabolic shock fronts develop a cusp from 
which modulated lump solutions emerge.   We have justified this behaviour with the observation that the Whitham modulation equations near the solitonic front are not hyperbolic.
The scaling of the maximum of 
the solution is linear with respect to the maximum amplitude of the 
initial data, and for the specific initial data considered, this scaling coefficient turns out to be about 10.
 Regarding the scaling of $x$ and $t$ as a function of $\epsilon$,  the  same scalings are observed as in the case of the  semiclassical limit of  focusing NLS. 

It would be interesting to identify the values of the break-up point 
$(x_{0},y_{0},t_{0})$ for given initial data. A way to obtain this 
information would be to solve the Whitham equations and to determine 
the point where their solutions develop a cusp for given initial 
data. A detailed study of the Whitham equations could also give an 
indication on how to make the above conjecture more precise, and how 
to prove it eventually.  Finally a more mathematical goal is the study of the integrability and Hamiltonian structure of the Whitham modulation equation as defined in \cite{Ferapontov}, \cite{Ferapontov1}.
This will be the subject of further work. 
\vskip 1cm
 \noindent
 {\bf\large Funding}
 T.G. acknowledges the support by the Leverhulme Trust Research Fellowship RF-2015-442.\\

{\bf Acknowledgements}
We thank Miguel Onorato, Peter Miller,  Karima Khusnutdinova for valuable discussions during the preparation of this manuscript.

\end{document}